\documentclass{article}
 \usepackage{graphicx}

\usepackage{amsmath}
\usepackage{amscd}
\usepackage{amstext}
\usepackage{amsbsy}
\usepackage{amsopn}
\usepackage{amsxtra}
\usepackage{amssymb}
\usepackage[latin1]{inputenc}
\usepackage{colortbl}

\newcommand{\Real}{\mathbb{R}}

\newcommand{\mysection}[1]{\vspace*{0cm}\section{#1}\vspace{0cm}}
\newcommand{\mysubsection}[1]{\vspace*{0cm}\subsection{#1}\vspace{0cm}}

\newcommand{\citeyearpar}[1]{\cite{#1}}
\newcommand{\citep}[1]{\cite{#1}}
\newcommand{\citet}[1]{\cite{#1}}

\newtheorem{lem}{Lemma}[section]

\newtheorem{defn}{Definition}[section]

\def\letters{0}
\def\Qhat{\widehat{Q}}
\def\Qhatsigma{\widehat{Q}}

  \renewenvironment{thebibliography}[1]{%
    \begin{oldthebibliography}{#1}%
      \setlength{\parskip}{0ex}%
      \setlength{\itemsep}{0ex}%
  }%
  {%
    \end{oldthebibliography}%
  }

\title{A quadratic measure of dependence}
\author{Sophie Achard}
\begin{document}

\maketitle

Brain Mapping Unit, University of Cambridge, Downing Site, Cambridge, CB2 3EB. email: sa428\makeatother cam.ac.uk\\

\vspace{1cm}

Abstract :\\

Asymptotic properties of a dimension-robust dependence measure are investigated. It is related to 
those used in independence tests, but is derivable, thus suitable for independent component analysis.
An adjustable kernel allows to accelerate the convergence of the estimator without affecting the 
bias.\\

Keywords :\\

dependence measure, U-statistics, characteristic functions, kernels.


\mysection{Introduction}

Since the 1950s, there has been a continuous research activity over the definition 
of measures of dependence, that is, positive functions that are equal to zero if and 
only if the variables are independent. The necessity of such measures first appeared 
in the construction of independence tests. Hoeffding \citeyearpar{hoeffding.1948.2} proposed 
to define an independence test by comparing the joint cumulative distribution 
function and the product of the marginal cumulative distribution functions. Then, 
in the 1970s several authors, including Rosenblatt \citeyearpar{rosenblatt.1975.1}, Blum {et al.} \citeyearpar{blum.1961.1},
Feuerverger \citeyearpar{feuerverger.1993.1}, have studied independence tests defined by 
comparing the joint density and the product of the marginal densities, or comparing 
the joint characteristic functions and the product of the marginal characteristic 
functions. But, in general, these tests are constructed to control the independence 
of only two variables, and are unsuitable in higher dimensions because of the curse 
of dimension in estimating the density. Recently, measures of dependence have received renewed interest as they play a 
crucial role in obtaining a procedure for independent component analysis (ICA) 
\citep{comon.1994.1, cardoso.1998.1}. This analysis aims at finding a transformation 
(usually linear) of a vector of observations, such that the transformed vector has 
independent components. To this end one minimizes a measure of the dependence 
between the transformed components. In order to employ efficient minimisation procedures, 
such a measure has to be differentiable with 
respect to 
the transformation, which is not the case for measures based on order statistics, 
for example. \\
In this letter, we study the dependence measure called quadratic dependence, introduced in 
\citet{achard.2003.1} and whose definition involves an adjustable kernel function. 
In order to relate this quadratic dependence with other 
existing dependence measures (e.g. the ones introduced by Chen et al. \citeyearpar{chen.2005.2}, 
Eriksson et al. \citeyearpar{eriksson.2002.2} and Kankainen \citeyearpar{kankainen.1995.1}), 
we derive two different expressions for it (Section \ref{sec.def}). The first one is 
based on the comparison 
of the joint characteristic function and the product of the marginal characteristic 
functions, which allows us to derive asymptotic properties of the estimator. 
The second one is based on its decomposition as U-statistics, and allows us to prove 
asymptotic normality and to gain insight on the crucial choice for the effect of the 
bandwidth of the kernel (Section \ref{sec.prop}). 
In Section \ref{sec.prop}, we apply the quadratic 
dependence in the context of independence tests.
\ifnum\letters=0
\section{Definitions of the quadractic dependence measure and estimations}\label{sec.def}
\fi
\ifnum\letters=1
\mysection{Definitions and estimations}\label{sec.def}
\fi
We introduce a dependence measure which is continuous and derivable, so as to allow convenient 
minimisation procedures. 
Let $\mathcal{K}$ be a summable function such that its Fourier transform is different 
from zero almost everywhere. Then, for any random variables $Y_1, \ldots, Y_K$, the equality of 
$E\left[\prod_{k=1}^{K}\mathcal{K}\left
  (y_k-Y_k\right )\right]$ and $\prod_{k=1}^{K}E\left[\mathcal{K}\left
  (y_k-Y_k\right )\right] $ for all vectors $(y_1,\ldots,y_K)$ in $\Real^K$ is equivalent to the independence of $Y_1, \ldots, Y_K$ \ifnum\letters=0 \citep{renyi.1966.2}\fi.
Thus, a dependence measure can be obtained by associating this characterization of dependence with a quadratic measure, as described below. This dependence measure is called the quadratic dependence and was first introduced by \citet{achard.2003.1}. 
\mysubsection{A kernel-based characterisation of independence}
\begin{defn}[Quadratic dependence]\label{def.DQ.kernel}
Let $\mathcal{K}$ be a square summable kernel function with Fourier
transform different from zero almost everywhere. For a set of $K$
random variables $Y_1,\ldots, Y_K$ (with finite variance), we define the quadratic measure of
their (mutual) dependence as
\[
Q(Y_1,\ldots, Y_K)=\frac{1}{2}\int D_{\mathbf{Y}}(y_1,\ldots, y_K)^2
 dy_1\ldots dy_K.
\]
where $\mathbf{Y} = (Y_1, \;\ \cdots ,Y_K)^T$ and
\begin{equation}\label{eq.depquadra.kernel}
D_{\mathbf{Y}}(y_1,\ldots, y_K) = E\left[\prod_{k=1}^{K}\mathcal{K}_h\left
  (y_k-\frac{Y_k}{\sigma_{Y_k}}\right )\right] -
\prod_{k=1}^{K}E\left[\mathcal{K}_h\left
  (y_k-\frac{Y_k}{\sigma_{Y_k}}\right )\right]
\end{equation}
where $\sigma_{Y_k}$ is a scale factor, that is, a positive functional of
the distribution of ${Y_k}$ such that $\sigma_{\lambda
Y_k}=|\lambda|\sigma_{Y_k}$, for all real constant $\lambda$, and $\mathcal{K}_h=\mathcal{K}(x/h)/h$.
\end{defn} 
\ifnum\letters=1 \vspace{-0.5cm} \fi
First of all, the following lemma establishes that the function $Q$ is a dependence measure. 
\begin{lem}\label{lem.equi.ind}
For any random variables $Y_1,\ldots, Y_K$, $Q(Y_1,\ldots, Y_K)=0$ if and only if the random variables $Y_1,\ldots, Y_K$ are independent.
\end{lem}
The proof of lemma \ref{lem.equi.ind} follows from the continuity of the characteristic functions and the equivalent expression of $Q$ in terms of the characteristic functions as stated in the lemma \ref{lem.equi.charac}.
\begin{lem}\label{lem.equi.charac}
Let us define $\psi_{\mathbf{Y}}$ the joint characteristic function of $Y_1, \ldots, Y_K$, $\psi_{Y_k}$ the characteristic function of $Y_k$, and $\psi_{\mathcal K}$ the Fourier
transform of $\mathcal K$. 
Let $D^c_{\mathbf{Y}}$ be the difference between the joint characteristic function and the product of the marginal characteristic functions :
\begin{equation}\label{psiformula}
D^c_{\mathbf{Y}}(y_1,\ldots, y_K) =
\psi_{\mathbf{Y}}(y_1,\ldots,y_K) - \prod_{k=1}^{K}\psi_{Y_k}(y_k)
\end{equation}
where $\mathbf{Y} = [Y_1, \;\ \cdots ,Y_K]^T$
Then, the quadratic dependence $Q$ can be expressed as a weighted average of $|D^c_{\mathbf{Y}}|^2$ :
\begin{eqnarray}
Q(Y_1,\ldots, Y_K) = \label{psidiffformula}
\frac{1}{2}\int
\prod_{k=1}^{K} \Big|
\frac{\sigma_{Y_k}\psi_{\mathcal{K}_h}(\sigma_{Y_k}y_k)}{2\pi}
      \Big|^2
|D^c_{\mathbf{Y}}(\mathbf{y})|^2 dy_1\ldots dy_K.
\end{eqnarray}
\end{lem}
\ifnum\letters=0
The proof of this lemma is given in appendix \ref{proof.lem.equi.charac}.
\fi
\ifnum\letters=1
This lemma is proved using the Parseval formula, that is, the Fourier transform is unitary. 
\fi
Also, it is easily verified from (\ref{eq.depquadra.kernel}) that the quadratic dependence is invariant by translation and by 
multiplication by a scalar.\\
The measure (\ref{psidiffformula}) has been considered by Kankainen \citeyearpar{kankainen.1995.1}, 
Eriksson et al.
\citeyearpar{eriksson.2002.2} and Feuerverger \citeyearpar{feuerverger.1993.1}, but only in the
particular case where $\mathcal{K}$ is a Gaussian kernel and without a scaling factor. It can also be seen as a generalisation of the measure defined by Rosenblatt 
\citeyearpar{rosenblatt.1975.1}. Indeed, when the bandwidth tends to zero and under 
usual hypotheses for the density and the kernel, $Q$ is equal to the 
quadratic measure of the difference between the joint density and the product of the 
marginal densities. \ifnum\letters=1 This result is stated in the following lemma and comes from Bochner's lemma and the Cauchy-Schwarz inequality.\fi
\begin{lem}\label{lem.depquadra.density}
Let us define $p_{\mathbf{Y}}$ the joint density of $Y_1,\ldots,Y_K$ and $p_{Y_k}$ the 
density of $Y_k$. Let us assume that the kernel $\mathcal{K}$ is a Parzen-Rosenblatt kernel, that is
 $\lim_{|x|\rightarrow \infty}|x|\mathcal{K}(x)=0$. Then, for all $\mathbf{y}$ where the joint density is continuous, 
$$\lim_{h\rightarrow 0} E\left [\prod_{k=1}^K\mathcal{K}_h(y_k-Y_k/\sigma_{Y_k})\right ]=p_{\mathbf{Y}}(y_1/\sigma_{Y_1},\ldots,y_K/\sigma_{Y_K})/\prod_{k=1}^K\sigma_k.$$ And for all  $y_k$ where the marginal density $p_{Y_k}$ is continuous, 
$$\lim_{h\rightarrow 0} E\left [\mathcal{K}_h(y_k-Y_k/\sigma_{Y_k})\right ]=p_{Y_k}(y_k/\sigma_{Y_k})/\sigma_{Y_k}$$ 
Moreover, if $p_{\mathbf{Y}}$ is continuously differentiable up to the second order 
with bounded derivatives and the first moment of $\mathcal{K}$ and $\mathcal{K}^2$ exist, then
$$\lim_{h\rightarrow 0} Q(Y_1,\ldots,Y_K)=1/2\int |p_{\mathbf{Y}}(\mathbf{y})-\prod_{k=1}^K p_{Y_k}(y_k)|^2d\mathbf{y}.$$ 
\end{lem}  
\ifnum\letters=0 This result is proved in appendix \ref{proof.lem.depquadra.density}.\\ \fi
Chen and Bickel \citeyearpar{chen.2005.2} studied the minimum of an estimator of 
(\ref{psidiffformula}) in the context of linear ICA and proved its consistency independently 
of the choice of a kernel\footnote{Note also that, in their definition, 
there is no scaling factor in the weight function, and therefore, no invariance by multiplication 
by a factor}. They point however that the choice of a kernel and especially, variations of its bandwidth, 
can change dramatically the variance and convergence in moment of the 
estimators.
One purpose of the present study is to shed some light on the influence of the 
bandwidth on the behaviour of the quadratic dependence measure in the context of 
independence tests. 
\mysubsection{Kernel trick}
The quadratic dependence as rewritten in (\ref{psidiffformula}) is not easy to estimate because of the multiple integration. The following lemma derives a formula for the quadratic 
dependence from which a convenient estimator arises. The trick employed for this is
 specific to this measure, and is a first step to address the problem of the curse 
of dimension. 
\begin{lem}\label{lem.kernel.trick}
Let $\mathcal{K}_2$ be the convolution of $\mathcal{K}$ with its mirror, i.e. $\mathcal{K}_2(u)=\int \mathcal{K}(u+v)\mathcal{K}(v) dv$. For a set of $K$
random variables $Y_1,\ldots, Y_K$ (with finite variance), the quadratic measure of
their (mutual) dependence is equivalent to
\begin{eqnarray}\label{formula.DQ.kernel}
Q(Y_1,\ldots,Y_K) & = & \frac{1}{2}\left
\{E \left [\pi_\mathbf{Y}(\mathbf{Y})\right ] +
\prod_{k=1}^{K}E\left [\pi_{Y_k}(Y_k)\right ] -
2 E\left [\prod_{k=1}^{K}\pi_{Y_k}(Y_k)\right ]\right \},
\end{eqnarray}
\ifnum\letters=0 where 
\begin{eqnarray*}
\pi_\mathbf{Y}(\mathbf{y}) & = &
E\left [\prod_{i=1}^{K}
\mathcal{K}_{2,h}\left(\frac{y_i-Y_i(n)}{\sigma_{Y_i}}\right)\right ],\\
\pi_{Y_k}(y_k)& = &
E\left [\mathcal{K}_{2,h}\left(\frac{y_k-Y_k(n)}{\sigma_{Y_k}}\right) \right ],
\end{eqnarray*}
\fi
\ifnum\letters=1 where 
\begin{equation*}
\pi_\mathbf{Y}(\mathbf{y})  = 
E\left [\prod_{i=1}^{K}
\mathcal{K}_{2,h}\left(\frac{y_i-Y_i(n)}{\sigma_{Y_i}}\right)\right ],\quad \pi_{Y_k}(y_k) = 
E\left [\mathcal{K}_{2,h}\left(\frac{y_k-Y_k(n)}{\sigma_{Y_k}}\right) \right ],
\end{equation*}
\fi
and $\sigma_{Y_k}$ is a scale factor (see definition \ref{def.DQ.kernel}).
\end{lem}
\ifnum\letters=0
The proof of this lemma is given in appendix \ref{proof.lem.kernel.trick}.
\fi
\ifnum\letters=1
The proof of this lemma is obtained by developping the square term and using 
Fubini's theorem to inverse the integration variables.
\fi
This lemma shows that $Q$
depends on $\mathcal K$ only indirectly through ${\mathcal K}_2$, therefore we
can choose ${\mathcal K}_2$ directly without ever considering $\mathcal K$.
For consistency with its definition, ${\mathcal
  K}_2$ must be choosen such that its Fourier transform is a
positive summable even function, since its Fourier transform
corresponds to $|\psi_{\mathcal K}|^2$ where $\psi_{\mathcal K}$ is the Fourier
transform of a real square summable function. Moreover, the Fourier transform of $\mathcal{K}_2$ 
has to be different from zero almost everywhere.\\
\begin{figure}[h]
\begin{minipage}{4cm}
\textbf{(a)}\\

\vspace*{-0cm}\includegraphics[width=\textwidth,angle=270]{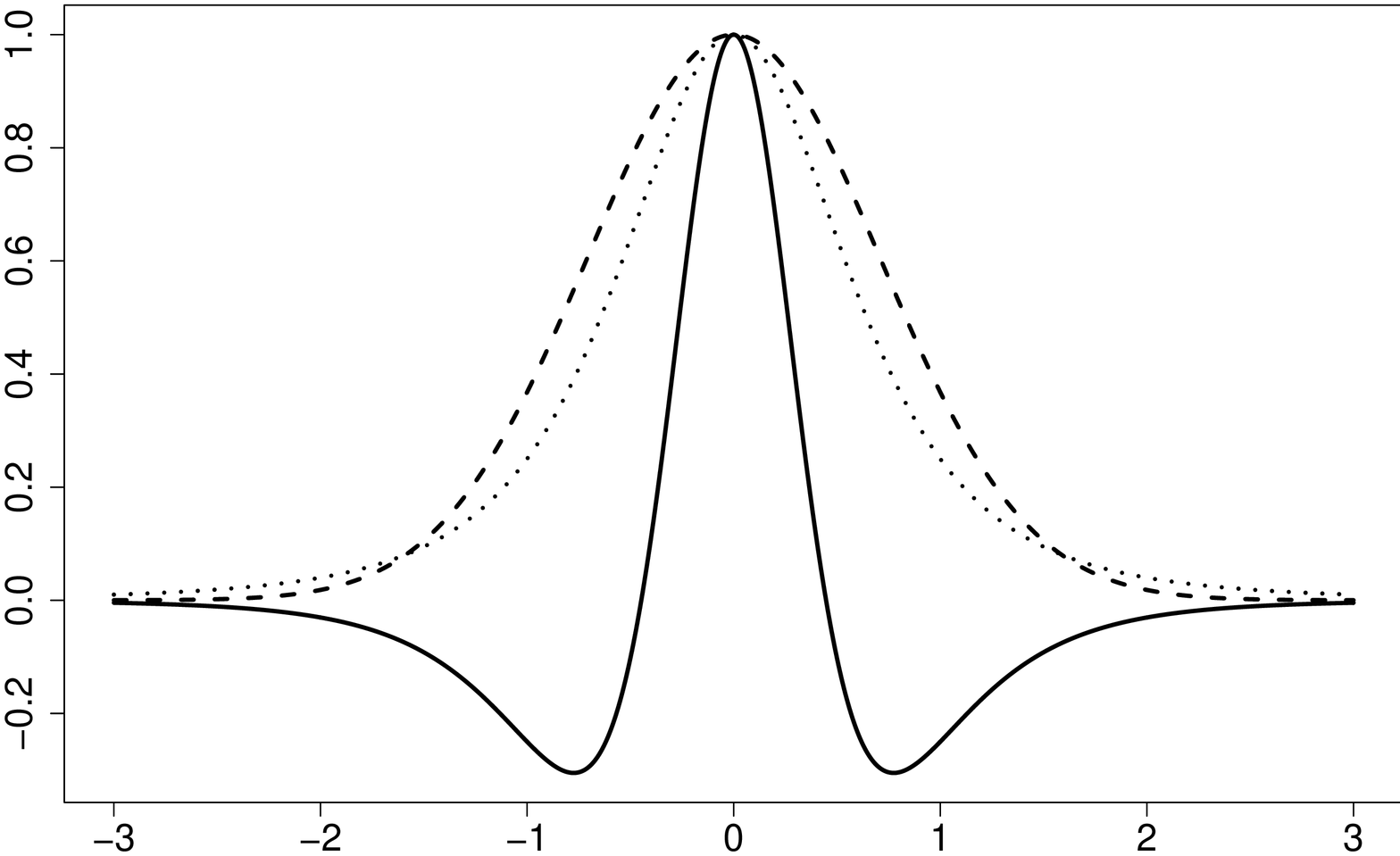}\end{minipage}  
\hspace*{3cm} \begin{minipage}{4cm}
\textbf{(b)}\\

\vspace*{-0cm}\includegraphics[width=\textwidth,angle=270]{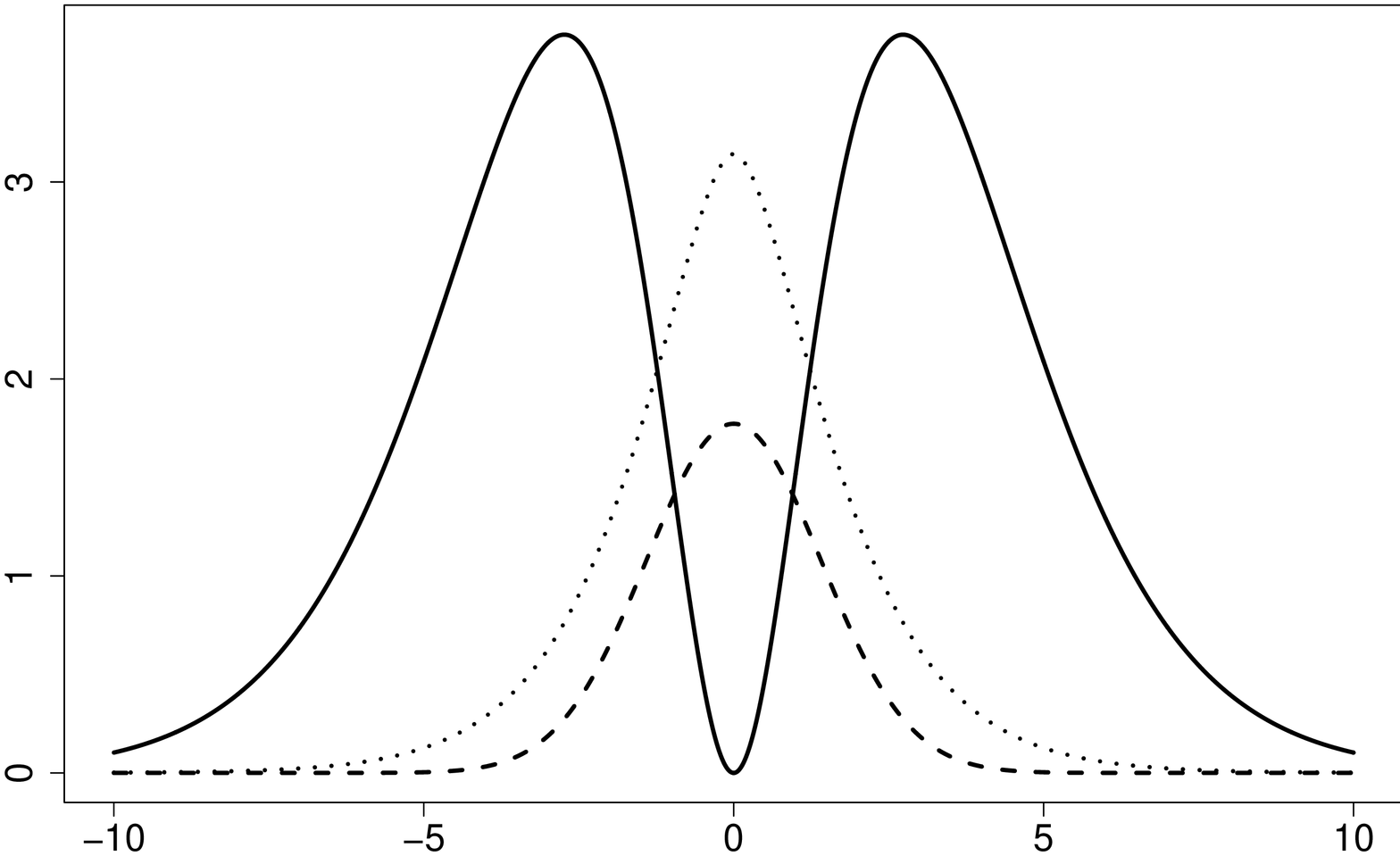}\end{minipage}\\
\caption{Plots of possible kernels and their Fourier transform, $- - -$, Gaussian. $\cdots$, square Cauchy. ---, the negative of the second derivative of the square Cauchy kernel. (a) is the plot of the kernels and (b) is the plot of their Fourier transform.}\label{fig.noyaux}
\end{figure}
Some possible choices for $\mathcal{K}_2$ are 
(denoting $\psi_{\mathcal{K}_2}$ the Fourier transform of $\mathcal{K}_2$) :
the Gaussian kernel, $\mathcal{K}_2(x) = e^{-x^2}$, $\psi_{\mathcal{K}_2}(t) = \sqrt{\pi}e^{-t^2/4}$, the square Cauchy kernel, $\mathcal{K}_2(x) = 1/(1+x^2)^2$, $\psi_{\mathcal{K}_2}(t) = \pi (|t|+1)e^{-|t|}$, the negative of the second derivative of the square Cauchy kernel, $\mathcal{K}_2(x) = -(20x^2-4)/(1+x^2)^4$, \medskip $\psi_{\mathcal{K}_2}(t) =4t^2\pi^3 (|t|+1)e^{-|t|}$.\\
The first two kernels correspond to density kernels (after normalizing) and differ only by their tail behaviour. The third kernel is not a density kernel and takes negative values (see figure \ref{fig.noyaux}). One may note that the kernel
$\mathcal{K}_2$ is related to Mercer kernels which are used especially
in Support Vector Machine\ifnum\letters=0 \citep{vapnik.1998.1}\fi.

\mysubsection{Estimation}
An estimator of $Q$ is defined using expression (\ref{formula.DQ.kernel}).  
In the sequel, the observed data will be denoted by $Y_k(n)$,
$n=1,\ldots,N$; $k=1, \ldots, K$; $N$ being the sample size and the scaling 
factor $\sigma=(\sigma_{Y_1},\ldots,\sigma_{Y_K})$ is supposed to be known, that is independent of the sample. 
Let us remark that (\ref{formula.DQ.kernel}) involves only the
expectation operator $E$, thus a natural estimator of $Q$ can be
obtained simply by replacing this operator with the sample average
$\widehat{E}$, defined as $\widehat{E}\phi(\mathbf{Y}) =
\sum_{n=1}^{N}\phi(\mathbf{Y}(n))/N$, for any function $\phi$ of 
 $K$ (real) variables. Thus, an estimator of $Q$ will be,
\begin{eqnarray}\label{formula.est.DQ}
\Qhat(Y_1,\ldots,Y_K)&= & \frac{1}{2}\left
\{\widehat{E}\widehat{\pi}_\mathbf{Y}(\mathbf{Y}) +
\prod_{k=1}^{K}\widehat{E}\widehat{\pi}_{Y_k}(Y_k) -
2\widehat{E}\prod_{k=1}^{K}\widehat{\pi}_{Y_k}(Y_k)\right \},
\end{eqnarray}
\ifnum\letters=0 where
\begin{eqnarray*}
\widehat{\pi}_\mathbf{Y}(\mathbf{y}) & = &
\frac{1}{N}\sum_{n=1}^{N}\prod_{i=1}^{K}
\mathcal{K}_{2,h}\left(\frac{y_i-Y_i(n)}{{\sigma}_{Y_i}}\right),\\
\widehat{\pi}_{Y_k}(y_k) & = &
\frac{1}{N}\sum_{n=1}^{N}
\mathcal{K}_{2,h}\left(\frac{y_k-Y_k(n)}{{\sigma}_{Y_k}}\right),
\end{eqnarray*}
\fi
\ifnum\letters=1 
where
\begin{equation*}
\widehat{\pi}_\mathbf{Y}(\mathbf{y})  = 
\frac{1}{N}\sum_{n=1}^{N}\prod_{i=1}^{K}
\mathcal{K}_{2,h}\left(\frac{y_i-Y_i(n)}{{\sigma}_{Y_i}}\right),\quad
\widehat{\pi}_{Y_k}(y_k)  = 
\frac{1}{N}\sum_{n=1}^{N}
\mathcal{K}_{2,h}\left(\frac{y_k-Y_k(n)}{{\sigma}_{Y_k}}\right),
\end{equation*}
\fi
Note that the computational cost of the estimator $\Qhat$ is 
of order $KN^2$.\\
As the exact expression of $Q$ is given in terms of the characteristic functions, the estimator $\Qhat$ can alternatively be rewritten in terms of the estimators of the characteristic functions.
\begin{lem}\label{lem.est.DQ.charac}
\begin{eqnarray}\label{formula.est.DQ.charac}
\hspace{-1.3cm}\Qhat(Y_1,\ldots, Y_K) & = & \frac{\prod_{k=1}^{K}{\sigma}_{Y_k}}{(2\pi)^K}\int \prod_{k=1}^{K}\psi_{\mathcal{K}_{2,h}}({\sigma}_{Y_k}y_k) \big |\widehat{\psi}_{Y}(\mathbf{y})- \prod_{k=1}^{K}\widehat{\psi}_{Y_k}(y_k)\big |^2 dy_1\ldots dy_K,
\end{eqnarray}
\ifnum\letters=0
\begin{eqnarray*}
\widehat{\psi}_{\mathbf{Y}}(y_1,\ldots,y_K) & = & \widehat{E} \left [\prod_{
k=1}^{K}e^{iy_kY_k}\right ] = \frac{1}{N}\sum_{n=1}^{N}\prod_{k=1}^{K}e^{iy_
kY_k(n)}\\
 \widehat{\psi}_{Y_k}(Y_k) & = & \widehat{E} \left [e^{iy_kY_k}\right ] = \frac{1}{N}\sum_{n=1}^{N}e^{iy_kY_k(n)}\\
\end{eqnarray*}
\fi
\ifnum\letters=1
$$ \widehat{\psi}_{\mathbf{Y}}(y_1,\ldots,y_K)  = \frac{1}{N}\sum_{n=1}^{N}\prod_{k=1}^{K}e^{iy_
kY_k(n)}, \quad \widehat{\psi}_{Y_k}(Y_k) = \frac{1}{N}\sum_{n=1}^{N}e^{iy_kY_k(n)}$$
\fi
\end{lem} 
\ifnum\letters=0
The proof is given in appendix \ref{proof.lem.est.DQ.charac}.
\fi
\ifnum\letters=1
This lemma is proved using the Parseval formula, as in lemma \ref{lem.equi.charac}. 
\fi

\mysection{Asymptotic properties}\label{sec.prop}
Having proposed a convenient estimator of the quadratic dependence, the objective is now to show 
asymptotic properties of the estimator $\Qhat$, in order to control its efficiency in independence tests.
First, we note that this estimator $\Qhatsigma$ 
can be expressed in terms of U-statistics. 
The asymptotic behaviour of the estimator $\Qhat$ under the 
hypothesis of dependence of the random variables is given first.
Then, using U-statistics, the variance of the estimator $\Qhatsigma$ is computed. 
Finally, it is shown that the estimator $\Qhatsigma$ converges to a Gaussian 
distribution.
\mysubsection{Convergence under the hypothesis of dependence}
\begin{lem}\label{lem.DQ-conv}
Suppose that the Fourier transform of $\mathcal{K}_2$ is positive, different from zero almost everywhere. 
Then, under the hypothesis of the dependence of the random variables $Y_1,\ldots, Y_K$,
$\lim_{N \rightarrow +\infty}\Qhat(Y_1,\ldots, Y_K)>0~~a.s.$,
for any cumulative distribution function of $\mathbf{Y}$.
\end{lem}
\ifnum\letters=0
The proof is given in appendix \ref{proof.lem.DQ-conv}.
\fi
\ifnum\letters=1
\vspace{-0.5cm}This lemma is proved using two properties of the estimator of the characteristic functions 
given by Csörg\H{o}, (1981),
$\sup_{\mathbf{y}\in B}\Big |
\widehat{\psi}_{\mathbf{Y}}(\mathbf{y})-\psi_{\mathbf{Y}}(\mathbf{y})\Big |
\xrightarrow{a.s.} 0, N \rightarrow + \infty$
and
$\sup_{\mathbf{y}\in B}\Big |
\prod_{k=1}^K\widehat{\psi}_{Y_k}(y_k)-\prod_{k=1}^K
\psi_{Y_k}(y_k)\Big |
\xrightarrow{a.s.} 0, N  \rightarrow + \infty$. 
\fi
\mysubsection{Bias and variance}
Unlike what is the case for the estimation of the density, the estimator $\Qhat$ is unbiased, 
that is $E[\Qhat]=Q$. This comes from Hoeffding \citeyearpar{hoeffding.1948.1}. This result is 
completely 
independent of the choice of the kernel.  Consequently, the bandwidth does not have to assume a 
specific dependence on the sample size in order to achieve convergence in mean. In particular, the 
bandwidth does not have to vanish as the sample size tends to infinity. Moreover, the 
convergence of the estimator of 
$\Qhat$ does not suffer from the problem of curse of dimensionality.\\
Also we show that the variance of the estimator $\Qhat$ goes to zero for any fixed bandwidth.
More precisely, the exact development of the variance is given below, following the 
development of Hoeffding \citeyearpar{hoeffding.1948.1}, we deduce the dominant 
terms in the expansion of the variance of $\Qhatsigma$, \ifnum\letters=0(the proof 
is given in appendix \ref{proof.lem.variance})\fi
\begin{eqnarray*}
\mathrm{var}(\Qhatsigma) & = & \frac{4}{N}\Sigma_{(11)}+\frac{4K^2}{N}\Sigma_{(22)}+4\frac{(K+1)^2}{N}\Sigma_{(33)}-4\frac{2(K+1)}{N}\Sigma_{(13)} \\
  &  & -4\frac{2K(K+1)}{N}\Sigma_{(23)}+2\frac{4K}{N}\Sigma_{(12)}+o(1/N)\\
\end{eqnarray*}
\ifnum\letters=0
$$\Sigma_{(11)}=E[\pi_{\mathbf{Y}}(\mathbf{Y})^2]-\theta_1^2$$
\begin{eqnarray*}
\Sigma_{(12)} & = & \Sigma_{(21)}=\frac{1}{K}\sum_{l=1}^KE[\prod_{k\neq l}E(\pi_{Y_k}(Y_k))\pi_{Y_l}(Y_l)\pi_{\mathbf{Y}}(\mathbf{Y})]-\theta_1\theta_2
\end{eqnarray*}
\begin{eqnarray*}
\Sigma_{(13)}=\Sigma_{(31)} & = & \frac{1}{K+1}E[\pi_{\mathbf{Y}}(\mathbf{Y})\prod_{k=1}^K\pi_{Y_k}(Y_k)]-\theta_1\theta_3\\
 & - & \frac{1}{K+1}\sum_{l=1}^KE\left
 [\pi_{\mathbf{Y}}(\mathbf{Y})\widetilde{\pi}_{Y_l}(Y_l)\right ]
\end{eqnarray*}
with $\widetilde{\pi}_{Y_l}(y_l)=E\left [\prod_{k\neq l}\pi_{Y_k}(Y_k)\mathcal{K}_{2,h}\left (\frac{Y_l-y_l}{\sigma_{Y_l}}\right )\right ]$ 
\begin{eqnarray*}
\Sigma_{(23)} & = & \Sigma_{(32)} = \frac{1}{K(K+1)}\sum_{k=1}^K\prod_{l
  \neq k}E[\pi_{Y_k}(Y_k)]E[\pi_{Y_l}(Y_l)\prod_{k=1}^K\pi_{Y_k}(Y_k)]-\theta_2\theta_3\\
& + & \frac{1}{K(K+1)}\sum_{l,m=1}^K\prod_{l
  \neq k}E[\pi_{Y_k}(Y_k)]E[\pi_{Y_l}(Y_l)\widetilde{\pi}_{Y_m}(Y_m)]\\
\
\end{eqnarray*}
\begin{eqnarray*}
\Sigma_{(22)} & = & \frac{1}{K^2}\sum_{l,m=1}^K\prod_{k\neq l}E[\pi_{Y_k}(Y_k)]\prod_{k\neq m}E[\pi_{Y_k}(Y_k)]E[\pi_{Y_l}(Y_l)\pi_{Y_m}(Y_m)\
]-\theta_2^2
\end{eqnarray*}
\begin{eqnarray*}
\Sigma_{(33)} & = & \frac{1}{(K+1)^2}E[\prod_{k=1}^K\pi^2_{Y_k}(Y_k)]+
\frac{1}{(K+1)^2}\sum_{l,m=1}^KE[\widetilde{\pi}_{Y_l}(Y_l)\widetilde{\
\pi}_{Y_m}(Y_m)]\\
& +&
\frac{2}{(K+1)^2}\sum_{l}^KE[\prod_{k=1}^K\pi_{Y_k}(Y_k)\widetilde{\pi\
}_{Y_l}(Y_l)]
- \theta_3^2\\
\end{eqnarray*}
\fi
\ifnum\letters=1
\vspace*{-1.7cm}
$$\mathrm{where,} \quad \Sigma_{(11)}=E[\pi_{\mathbf{Y}}(\mathbf{Y})^2]-\theta_1^2,\quad
\Sigma_{(12)}=\frac{1}{K}\sum_{l=1}^KE[\prod_{k\neq l}E(\pi_{Y_k}(Y_k))\pi_{Y_l}(Y_l)\pi_{\mathbf{Y}}(\mathbf{Y})]-\theta_1\theta_2$$
$$\Sigma_{(13)} =  \frac{1}{K+1}E[\pi_{\mathbf{Y}}(\mathbf{Y})\prod_{k=1}^K\pi_{Y_k}(Y_k)]-\theta_1\theta_3 - \frac{1}{K+1}\sum_{l=1}^KE\left
 [\pi_{\mathbf{Y}}(\mathbf{Y})\widetilde{\pi}_{Y_l}(Y_l)\right ]$$

\vspace{-0.5cm}
\begin{eqnarray*}
\Sigma_{(23)} & =  & \frac{1}{K(K+1)}\sum_{k=1}^K\prod_{l
  \neq k}E[\pi_{Y_k}(Y_k)]E[\pi_{Y_l}(Y_l)\prod_{k=1}^K\pi_{Y_k}(Y_k)]-\theta_2\theta_3\\
& + & \frac{1}{K(K+1)}\sum_{l,m=1}^K\prod_{l
  \neq k}E[\pi_{Y_k}(Y_k)]E[\pi_{Y_l}(Y_l)\widetilde{\pi}_{Y_m}(Y_m)]\\
\end{eqnarray*}
$$
\Sigma_{(22)}  =  \frac{1}{K^2}\sum_{l,m=1}^K\prod_{k\neq l}E[\pi_{Y_k}(Y_k)]\prod_{k\neq m}E[\pi_{Y_k}(Y_k)]E[\pi_{Y_l}(Y_l)\pi_{Y_m}(Y_m)\
]-\theta_2^2$$
\begin{eqnarray*}
\Sigma_{(33)} & = & \frac{1}{(K+1)^2}E[\prod_{k=1}^K\pi^2_{Y_k}(Y_k)]+
\frac{1}{(K+1)^2}\sum_{l,m=1}^KE[\widetilde{\pi}_{Y_l}(Y_l)\widetilde{\
\pi}_{Y_m}(Y_m)]\\
& + &
\frac{2}{(K+1)^2}\sum_{l}^KE[\prod_{k=1}^K\pi_{Y_k}(Y_k)\widetilde{\pi\
}_{Y_l}(Y_l)]
- \theta_3^2\\
\end{eqnarray*}
with $\widetilde{\pi}_{Y_l}(y_l)=E\left [\prod_{k\neq l}\pi_{Y_k}(Y_k)\mathcal{K}_{2,h}\left (\frac{Y_l-y_l}{\sigma_{Y_l}}\right )\right ]$.\\
\fi
None of these quantities depend on $N$, but they are dependent on the 
choice of the kernel and its bandwidth, and on the distribution of the observations.

\mysubsection{Asymptotic Gaussian distribution and hypothesis test}

\begin{figure}[h]
\includegraphics[width=0.5\textwidth,angle=270]{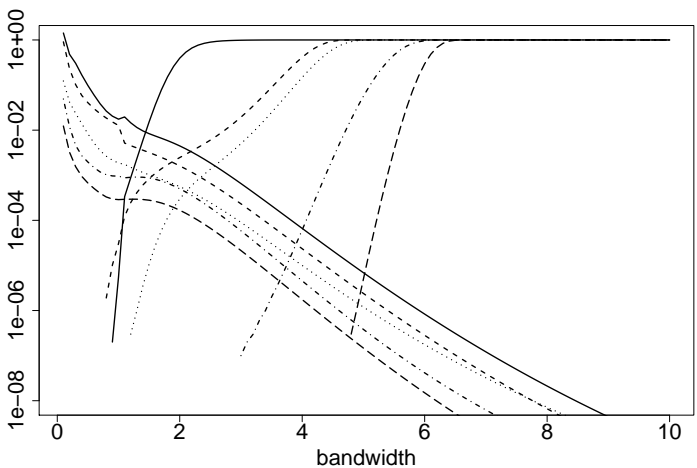}
\caption{Variance (decreasing) and type II error (increasing) for different sample sizes. ---, $N=100$, - - -, $N=200$, $\cdots$, $N=400$, - $\cdot$ -, $N=800$, -- -- --, $N=1600$. }\label{fig.bandwidth}
\end{figure}

The quadratic dependence measure studied in this paper provides us with an estimator for 
the evaluation of 
the dependence between variables. In the following, we construct a hypothesis 
test of independence based on the quadratic dependence measure. The asymptotic laws under 
the hypotheses of independence (denoted H$_0$) and dependence (denoted H$_1$) are deduced.
Finally, it is shown that this hypothesis test of independence is consistent for 
any choice of the bandwidth.\\
$\bullet$ \underline{Law under the hypothesis of independence} (denoted H$_0$):\\
The estimator
$N\Qhatsigma$ follows asymptotically a law of $\gamma\chi^2(\beta)$ where $\gamma$ and
$\beta$ are $\gamma=V_1/2E_1 ~ \mathrm{and} ~ \beta= 2E_1^2/V_1$, where $E_1=\lim_{N \rightarrow \infty}NE[\Qhatsigma]$ under H$_0$, and $V_1=\lim_{N \rightarrow \infty} N\mathrm{var}(\Qhatsigma)$ 
under H$_0$. It holds that 
$$E_1=\prod_{k=1}^{K}\int \mathcal{K}_{2,h}^2(x)dx -
\prod_{k=1}^{K}E[\pi_{Y_k}(Y_k)]-\sum_{k=1}^K(\int
\mathcal{K}_{2,h}^2(x)dx-E[\pi_{Y_k}(Y_k)])\prod_{l=1,l\neq k}^{K}E[\pi_{Y_l\
}(Y_l)]
$$
and
\begin{eqnarray*}
V_1 & = &
2\prod_{k=1}^{K}E[\pi_{Y_k}(Y_k)]^2-4\prod_{k=1}^{K}E[\pi_{Y_k}(Y_k)^2\
]+4\prod_{k=1}^{K}E[\pi_{Y_k}(Y_k)]\\
 & + &
2\sum_{k=1}^K(E[\pi_{Y_k}(Y_k)]-E[\pi_{Y_k}(Y_k)]^2)\prod_{l=1,l\neq
  k}^{K}E[\pi_{Y_l}(Y_l)]^2\\
 &  - & 4\sum_{k=1}^K(E[\pi_{Y_k}(Y_k)]-E[\pi_{Y_k}(Y_k)^2])\prod_{l\
=1,l\neq k}^{K}E[\pi_{Y_l}(Y_l)^2]\\
 & + & 2 \sum_{k=1}^K \sum_{m=1,m\neq
  k}^K(E[\pi_{Y_k}(Y_k)]^2E[\pi_{Y_m}(Y_m)]^2-2E[\pi_{Y_k}(Y_k)^2]E[\pi_{Y_m}(Y_m)]^2\\
 & + & E[\pi_{Y_k}(Y_k)^2]E[\pi_{Y_m}(Y_m)^2])\prod_{l=1,l\neq k,m}^{K\
}E[\pi_{Y_l}(Y_l)]
\end{eqnarray*}
This result is due to Kankainen \citeyearpar{kankainen.1995.1}.\\ 
\noindent $\bullet$ \underline{Law under the hypothesis of dependence} (denoted H$_1$):\\
$\sqrt{N}(\Qhatsigma-Q)$ follows asymptotically a normal law with zero mean and $\widetilde{\Sigma}$ variance,
where $\widetilde{\Sigma}$ is,
$$\hspace*{0cm}
\widetilde{\Sigma}  =  4\Sigma_{(11)}+4K^2\Sigma_{(22)}+4(K+1)^2\Sigma_{(33)}-8(K+1)\Sigma_{(13)}-8K(K+1)\Sigma_{(23)}+8K\Sigma_{(12)}$$
with $\Sigma$ the variance-covariance matrix of the corresponding U-statistics, which 
depends on $\mathcal{K}_2$ and $h$. 
\begin{lem}\label{lem.ind}
The independence test defined above is consistent for any choice of the bandwidth :
Given $\alpha$, the level of significance, we define $q_{\alpha}$ the smallest number satisfying the inequality $P_{\mathrm{H}_0}(\Qhat>q_{\alpha})=1-F_{\gamma\chi^2(\beta)}(Nq_{\alpha}) \leq \alpha$.\\
Then, the power of the test $1-P_{\mathrm{H}_1}(\Qhat<q_{\alpha})$ tends to 1 as $N$ goes to infinity. \\
In addition, the power of the test admits a lower bound :
\begin{equation}\label{eq.power.bound}
1-P_{\mathrm{H}_1}(\Qhat<q_{\alpha})=P_{\mathrm{H}_1}(\Qhat>q_{\alpha})> 1-\frac{ \mathrm{var}(\Qhat)}{q_\alpha-Q}
\end{equation}
\end{lem} 
\ifnum\letters=0
The proof is given in appendix \ref{proof.lem.ind}. Note that the lower bound in 
(\ref{eq.power.bound}) is not sharp as is illustrated in figure \ref{fig.bandwidth}.
\fi
\ifnum\letters=1
This lemma is proved using the Chebychev inequality and the asymptotic properties of the variance. 
Note that the lower bound in 
(\ref{eq.power.bound}) is not sharp, as is illustrated in figure \ref{fig.bandwidth}.
\fi
\mysubsection{Convergence rate and choice of the bandwidth}
As the asymptotic bias of the estimator tends to zero when the size of the sample 
$N$ goes to infinity without any constraint on the bandwidth, it is not necessary to 
choose the bandwidth so as to make 
a tradeoff between the bias and the variance. As a result, the bandwidth can rather 
be adjusted in a tradeoff between the minimisation of the variance and the minimisation of the 
type II error of the test.\\
With the expression of the variance given above, it is clear that when the bandwidth $h$ of 
the kernel increases, the variance of the estimator will decrease. But, as the bandwidth $h$ 
increases, the type II error of the test is expected to increase. Indeed, the asymptotic power 
of the test is defined by $$1-P_{\mathrm{H}_1}(\Qhat<q_{\alpha})=1-\Phi \left(\frac{(q_{\alpha}-Q)\sqrt{N}}{\sigma}\right )$$ where for a given $\alpha$, $q_{\alpha}$ verifies 
$P_{\mathrm{H}_0}(\Qhat>q_{\alpha})=1-F_{\gamma\chi^2(\beta)}(Nq_{\alpha}).$
Figure \ref{fig.bandwidth} illustrates the behaviour of optimal choices of the 
bandwidth depending on the size of the sample. On this figure, we observe that for 
a given bandwidth the convergence of the variance of the estimator is slow (of order $1/N$). 
But, if the bandwidth is adjusted to make a tradeoff between the variance and the power of 
the test, the convergence rate is tremendously increased. Future work has to be done to 
quantify this increase and to propose a computational rule to optimize the bandwidth.

\section*{Acknowledgements}
The author would like to thank D.-T. Pham for pointing out some theoretical problems.

\mysection{Conclusion}
The quadratic dependence measure is revisited and its asymptotic properties 
are demonstrated. The convergence rate in terms of the variance is of order 
$1/N$, and the power of the test defined by this measure converges to one with 
a rate of $1/N$ at least, $N$ being the sample size. \\
The introduction of a kernel frame in the definition of the quadratic dependence 
measure enables us to propose an efficient estimator of computational cost of order 
$KN^2$, with $K$ the dimension of the problem. This kernel is adjusted with a bandwidth whose choice does not affect the bias, which differs with the case of the estimation of density. 
Because of this property, the bandwidth can be chosen in terms of the sample size $N$, 
so as to increase the convergence rate in terms of the variance and of the power of the test rather 
than debiasing the estimator.

\section*{Appendix}

\appendix

\section{Proof of lemma \ref{lem.equi.charac}}\label{proof.lem.equi.charac}

After a change in the integral variable, $Q$ is written as, 
\[
\hspace{-1cm}
Q(Y_1,\ldots, Y_K)=\frac{1}{\prod_{k=1}^{K}\sigma_{Y_k}}\int \left
  \{E\left [\prod_{k=1}^{K}\mathcal{K}_h\left
  (\frac{y_k-Y_k}{\sigma_{Y_k}}\right )\right] -
\prod_{k=1}^{K}E\left[\mathcal{K}_h\left
  (\frac{y_k-Y_k}{\sigma_{Y_k}}\right )\right]\right \}^2
 dy_1\ldots dy_K.
\]

Let us denote a new kernel,
$$\mathcal{K}_{h,\sigma_{Y_k}}(x)=\mathcal{K}_h\left (\frac{x}{\sigma_{Y_k}}\right
)$$

then, $Q$ is expressed as,
\[
\hspace{-1cm}
Q(Y_1,\ldots, Y_K)=\frac{1}{\prod_{k=1}^{K}\sigma_{Y_k}}\int \left \{
  \left [\prod_{k=1}^{K}\mathcal{K}_{h,\sigma_{Y_k}}\right ]*dF_{\mathbf{Y}}(y_1,\ldots,y_K) -
\prod_{k=1}^{K} \mathcal{K}_{h,\sigma_{Y_k}}*dF_{Y_k}(y_k)\right \}^2
 d\mathbf{y}.
\]
where $d\mathbf{y}:=dy_1\ldots dy_K$.
Moreover, the Fourier transform of
$(y_1,\ldots,y_K) \mapsto \left [\prod_{k=1}^{K}\mathcal{K}_{h,\sigma_{Y_k}}\right
]*dF_{\mathbf{Y}}(y_1,\ldots,y_K)$ and $y_k \mapsto
\mathcal{K}_{h,\sigma_{Y_k}}*dF_{Y_k}(y_k)$ are respectively equal to 
\[
(y_1,\ldots,y_K) \mapsto \prod_{k=1}^K
\sigma_{Y_k}\psi_{\mathcal{K}_h}(\sigma_{Y_k}y_k)\psi_{\mathbf{Y}}(y_1,\ldots,y_K)
\text{ and } y_k \mapsto \sigma_{Y_k}\psi_{\mathcal{K}_h}(\sigma_{Y_k}y_k)\psi_{Y_k}(y_k\
)
\]
Then using the Parseval's formula, that is the Fourier transform is unitary, the lemma is proved.

\section{Proof of lemma \ref{lem.depquadra.density}}\label{proof.lem.depquadra.density}

The lemma \ref{lem.depquadra.density} comes directly from Bochner's lemma applied to the 
convolution of the kernel and the density function. The limit behaviour 
of $Q$ when $h$ tends to zero is proved by applying a Taylor development to the density.
The proof is rather computational, and we will only give the proof for,
$$\lim_{h\rightarrow 0}\left[\int \left [\left (\int \mathcal{K}(v)p_{Y_k}(hv+t)dv \right )^2-(p_{Y_k}(t))^2\right ]dt \right ]=0.$$

Indeed, by applying the triangular inequality and the Cauchy-Schwarz inequality, the desired 
quantity is bounded by, $hM_{p''_{Y_k}}\int |v|\mathcal{K}^2(v)dv$, which proves the limit.
The other terms behave similarly.

\section{Proof of lemma \ref{lem.kernel.trick}}\label{proof.lem.kernel.trick}

This lemma is proved by developping the square term in the definition of $Q$ see definition \ref{def.DQ.kernel}.
\begin{eqnarray*}
\lefteqn{D_{\mathbf{Y}}(y_1,\ldots, y_K)^2=}\\
 & &  \left \{ E\left[\prod_{k=1}^{K}\mathcal{K}_h\left
  (y_k-\frac{Y_k}{\sigma_{Y_k}}\right )\right] -
\prod_{k=1}^{K}E\left[\mathcal{K}_h\left
  (y_k-\frac{Y_k}{\sigma_{Y_k}}\right )\right]\right \}^2\\
 & = & \left \{ \int \prod_{k=1}^{K}\mathcal{K}_h\left
  (y_k-\frac{u_k}{\sigma_{Y_k}}\right
  )dF_{\mathbf{Y}}(u_1,\ldots,u_K)
  - \prod_{k=1}^{K}\int \mathcal{K}_h\left
  (y_k-\frac{u_k}{\sigma_{Y_k}}\right )dF_{Y_k}(u_k)\right \}^2\\
 & = & \iint \prod_{k=1}^{K}\mathcal{K}_h\left
  (y_k-\frac{u_k}{\sigma_{Y_k}}\right
  )\mathcal{K}_h\left
  (y_k-\frac{v_k}{\sigma_{Y_k}}\right
  )dF_{\mathbf{Y}}(\mathbf{u})dF_{\mathbf{Y}}(\mathbf{v})\\
 & + &
 \prod_{k=1}^{K}\iint \mathcal{K}_h\left
  (y_k-\frac{u_k}{\sigma_{Y_k}}\right )\mathcal{K}_h\left
  (y_k-\frac{v_k}{\sigma_{Y_k}}\right
  )dF_{Y_k}(u_k)dF_{Y_k}(v_k)\\
& - & 2 \int \prod_{k=1}^{K} \int \mathcal{K}_h\left
  (y_k-\frac{u_k}{\sigma_{Y_k}}\right )\mathcal{K}_h\left
  (y_k-\frac{v_k}{\sigma_{Y_k}}\right )  dF_{Y_k}(v_k)dF_{\mathbf{Y}}(\mathbf{u})\\
\end{eqnarray*}

So as to apply the Fubini's theorem, the properties 
of convergence of integrals have to be checked :

\begin{eqnarray*}
\lefteqn{\iiint \prod_{k=1}^{K}\left |\mathcal{K}_h\left
  (y_k-\frac{u_k}{\sigma_{Y_k}}\right
  )\mathcal{K}_h\left
  (y_k-\frac{v_k}{\sigma_{Y_k}}\right
  )\right
  |dF_{\mathbf{Y}}(\mathbf{u})dF_{\mathbf{Y}}(\mathbf{v})d\mathbf{y}}\\
 & \leq &  \iint \int \prod_{k=1}^{K}\left |\mathcal{K}_h\left
  (y_k-\frac{u_k}{\sigma_{Y_k}}\right
 )\mathcal{K}_h\left
  (y_k-\frac{v_k}{\sigma_{Y_k}}\right
  )\right |d\mathbf{y}dF_{\mathbf{y}}(\mathbf{u})dF_{\mathbf{Y}}(\mathbf{v})\\
& \leq & \iint dF_{\mathbf{Y}}(\mathbf{u})dF_{\mathbf{Y}}(\mathbf{v})
  \prod_{k=1}^{K} \int |\mathcal{K}_h|^2 < \infty\\
\end{eqnarray*}

then, for all $k=1,\ldots,K$,

\begin{eqnarray*}
\int \iint \left | \mathcal{K}_h\left
  (y_k-\frac{u_k}{\sigma_{Y_k}}\right )\mathcal{K}_h\left
  (y_k-\frac{v_k}{\sigma_{Y_k}}\right
  )\right |dF_{Y_k}(u_k)dF_{Y_k}(v_k)d\mathbf{y} & < & \infty \\
\end{eqnarray*}

and
\begin{eqnarray*}
\int \int \prod_{k=1}^{K} \int \left |\mathcal{K}_h\left
  (y_k-\frac{u_k}{\sigma_{Y_k}}\right )\mathcal{K}_h\left
  (y_k-\frac{v_k}{\sigma_{Y_k}}\right )\right |
  dF_{Y_k}(v_k)dF_{\mathbf{Y}}(\mathbf{u})d\mathbf{Y} & < & \infty\\
\end{eqnarray*}

This concludes the proof of the lemma by applying Fubini's theorem.

\section{Proof of lemma \ref{lem.est.DQ.charac}}\label{proof.lem.est.DQ.charac}

This lemma is proved in the same way as lemma \ref{lem.equi.charac}, by using the 
Parseval's formula and the following equality :

\begin{eqnarray*}
\lefteqn{\Qhat(Y_1,\ldots, Y_K) = \frac{\prod_{k=1}^{K}{\sigma}_{Y_k}}{(2\pi)^K}\int \prod_{k=1}^{K}\psi_{\mathcal{K}_{2,h}}({\sigma}_{Y_k}y_k) \big |\widehat{\psi}_{Y}(\mathbf{y})- \prod_{k=1}^{K}\widehat{\psi}_{Y_k}(y_k)\big |^2 dy_1\ldots dy_K}\\
 & = &
\frac{1}{N^2}\sum_{m=1}^N\sum_{n=1}^N\prod_{k=1}^K\int e^{iy_k\left
  (\frac{Y_k(n)-Y_k(m)}{{\sigma_{Y_k}}}\right )}\psi_{\mathcal{K}_{2,h}}(y_k)dy_k\\
& & +\prod_{k=1}^K\frac{1}{N^2}\sum_{m=1}^N\sum_{n=1}^N\int e^{iy_k\left
  (\frac{Y_k(n)-Y_k(m)}{{\sigma_{Y_k}}}\right )}\psi_{\mathcal{K}_{2,h}}(y_k)dy_k\\
& & -2\frac{1}{N}\sum_{m=1}^N\prod_{k=1}^K\frac{1}{N}\sum_{n=1}^N\int e^{iy_k\left
  (\frac{Y_k(n)-Y_k(m)}{{\sigma_{Y_k}}}\right )}\psi_{\mathcal{K}_{2,h}}(y_k)dy_k\\
\end{eqnarray*}

\section{Proof of lemma \ref{lem.DQ-conv}}\label{proof.lem.DQ-conv}

As the variables $Y_1,\ldots, Y_K$ are not independent, there exists $y_1,\ldots,y_K$ such that the following inequality is true : $$\psi_{\mathbf{Y}}(y_1,\ldots, y_K)\neq
\prod_{k=1}^K\psi_{Y_k}(y_k)$$

Then there exists a bounded open $U$ having positive Lebesgue measure such that 

$$\inf_{\mathbf{y} \in U}|\psi_{\mathbf{Y}}(y_1,\ldots, y_K)-
\prod_{k=1}^K\psi_{Y_k}(y_k)|>0$$

As the Fourier transform of $\mathcal{K}_2$ is different form zero allmost everywhere, we obtain the following inequality :

$$\int_{U} |\psi_{\mathbf{Y}}(y_1,\ldots, y_K)-
\prod_{k=1}^K\psi_{Y_k}(y_k)|^2\prod_{k=1}^K\psi_{\mathcal{K}_{2,h}}(y_k)d\mathbf{y}>0$$

From Csörg\H{o} \citeyearpar{csorgo.1981.2},

$$\sup_{\mathbf{y}\in B}\Big |
\widehat{\psi}_{\mathbf{Y}}(\mathbf{y})-\psi_{\mathbf{Y}}(\mathbf{y})\Big |
\xrightarrow{a.s.} 0, N \rightarrow + \infty$$
and
$$\sup_{\mathbf{y}\in B}\Big |
\prod_{k=1}^K\widehat{\psi}_{Y_k}(y_k)-\prod_{k=1}^K
\psi_{Y_k}(y_k)\Big |
\xrightarrow{a.s.} 0, N  \rightarrow + \infty$$

for any bounded set $B$, and in particular for $B=U$.

Thus, as the Fourier transforms of the kernel are summable,

$$\lim_{N \rightarrow +\infty}\int_{U} \Big |
\widehat{\psi}_{\mathbf{Y}}(\mathbf{y})-\psi_{\mathbf{Y}}(\mathbf{y})\Big
|\prod_{k=1}^K\psi_{\mathcal{K}_{2,h}}(y_k)d\mathbf{y}=0~~a.s.$$
and
$$\lim_{N \rightarrow +\infty}\int_{U} \Big |
\prod_{k=1}^K\widehat{\psi}_{Y_k}(y_k)-\prod_{k=1}^K
\psi_{Y_k}(y_k)\Big | \prod_{k=1}^K\psi_{\mathcal{K}_{2,h}}(y_k)d\mathbf{y}=0~~a.s.$$

Let us now remark that,

\begin{eqnarray*}
\lefteqn{\int_{U} |\widehat{\psi}_{\mathbf{Y}}(y_1,\ldots, y_K)-
\prod_{k=1}^K\widehat{\psi}_{Y_k}(y_k)|\prod_{k=1}^K\psi_{\mathcal{K}_{2,h}}(y_k)d\mathbf{y}\geq }\\
 & &\hspace*{0.5cm}\int_{U} |\psi_{\mathbf{Y}}(y_1,\ldots, y_K)-
\prod_{k=1}^K\psi_{Y_k}(y_k)|\prod_{k=1}^K\psi_{\mathcal{K}_{2,h}}(y_k)d\mathbf{y}\\
 & & -\int_{U} \Big |
\widehat{\psi}_{\mathbf{Y}}(\mathbf{y})-\psi_{\mathbf{Y}}(\mathbf{y})\Big |
\prod_{k=1}^K\psi_{\mathcal{K}_{2,h}}(y_k)d\mathbf{y}\\
 & - &\hspace*{0.5cm}\int_{U} \Big |
\prod_{k=1}^K\widehat{\psi}_{Y_k}(y_k)-\prod_{k=1}^K
\psi_{Y_k}(y_k)\Big |\prod_{k=1}^K\psi_{\mathcal{K}_{2,h}}(y_k)d\mathbf{y}\\
\end{eqnarray*}

which leads us to :

$$\liminf_{N \rightarrow +\infty}\int_{U} |\widehat{\psi}_{\mathbf{Y}}(y_1,\ldots, y_K)-
\prod_{k=1}^K\widehat{\psi}_{Y_k}(y_k)|\prod_{k=1}^K\psi_{\mathcal{K}_{2,h}}(y_k)d\mathbf{y}>0$$

With the Cauchy-Schwarz inequality, we get :

$$\liminf_{N \rightarrow +\infty}\int_{U} |\widehat{\psi}_{\mathbf{Y}}(y_1,\ldots, y_K)-
\prod_{k=1}^K\widehat{\psi}_{Y_k}(y_k)|^2\prod_{k=1}^{K}\psi_{\mathcal{K}_{2,h}}(y_k)d\mathbf{y}>0$$

Indeed,

\begin{eqnarray*}
\lefteqn{\int_{U} |\widehat{\psi}_{\mathbf{Y}}(y_1,\ldots, y_K)-
\prod_{k=1}^K\widehat{\psi}_{Y_k}(y_k)|^2\prod_{k=1}^{K}\psi_{\mathcal{K}_{2,h}}(y_k)d\mathbf{y}}\\
& \geq  & \frac{\int_{U} |\widehat{\psi}_{\mathbf{Y}}(y_1,\ldots, y_K)-
\prod_{k=1}^K\widehat{\psi}_{Y_k}(y_k)|\prod_{k=1}^K\psi_{\mathcal{K}_{2,h}}(y_k)d\mathbf{y}}{\int_{U} \prod_{k=1}^K\psi_{\mathcal{K}_{2,h}}(y_k)d\mathbf{y}}
\end{eqnarray*}

\section{U-statistics decomposition}\label{proof.U-stat}

Using the concept of U-statistics, $\Qhatsigma$ is expressed in a different way:

$\Qhatsigma=U'_1 + U'_2 -2U'_3$ where 

\begin{eqnarray*}
U'_1 & = & \frac{N-1}{N}U_1+\frac{b_N^{(1)}}{\sqrt{N}}\\
U'_2 & = & \frac{C_N^{2K}(2K)!}{N^{2K}}U_2+\frac{b_N^{(2)}}{\sqrt{N}}\\
U'_2 & = & \frac{C_N^{K+1}(K+1)!}{N^{K+1}}U_3+\frac{b_N^{(3)}}{\sqrt{N}}\\
\end{eqnarray*}

$b_N^{(1)}$, $b_N^{(2)}$, $b_N^{(3)}$ are random variables such that the limit of $E[b_N^{(i)}]^2$ is equal to zero if $N$ tends to infinity.

And $U_1$, $U_2$ and $U_3$ are the U-statistics associated to $U'_1$, $U'_2$ and $U'_3$ (the exact formulas are given at the end of this part).

In the following lemmas, we deduce the exact formula of $U_1$, $U_2$, $U_3$ and $b_N^{(1)}$, $b_N^{(2)}$ and $b_N^{(3)}$:

\begin{lem}
Let us decompose the estimation of $E \left [\pi_\mathbf{Y}(\mathbf{Y})\right ]$ in terms of U-statistics:
 
\begin{eqnarray*}
U'_1(\mathbf{Y}(1),\ldots,\mathbf{Y}(N)):=\widehat{E}[\widehat{\pi}_{\mathbf{Y}}(\mathbf{Y})] & = & \frac{1}{N}\sum_{n=1}^{N}\frac{1}{N}\sum_{m=1}^{N}\prod_{k=1}^{K}\mathcal{K}_{2,h}\left
 ( \frac{Y_k(n)-Y_k(m)}{\sigma_{Y_k}}\right )\\
 & = & \frac{N-1}{N}U_1+\frac{b_N^{(1)}}{\sqrt{N}} \\
\end{eqnarray*}
where 
$$U_1(\mathbf{Y}(1),\ldots,\mathbf{Y}(N)) =\frac{2}{N(N-1)}\sum_{1\leq i<j\leq N}
\prod_{k=1}^{K}\mathcal{K}_{2,h}\left
 ( \frac{Y_k(i)-Y_k(j)}{\sigma_{Y_k}}\right
 )$$
and 
$$b_N^{(1)}=\frac{1}{\sqrt{N}} \prod_{k=1}^{K}\mathcal{K}_{2,h}(0)$$
\end{lem}

\begin{lem}
Let us define the set $S_2=\{(i_1,\ldots,i_{K},j_1,\ldots,j_{K})| \forall k, 1\leq k \leq 2K, 1 \leq i_k \leq N,  1 \leq j_k \leq N\}$.
The second term in the expression of $\Qhatsigma$ corresponds to the estimation of $\prod_{k=1}^{K}E\left [\pi_{Y_k}(Y_k)\right ]$, and is written in terms of U-statistics as:

\begin{eqnarray*}
U'_2(\mathbf{Y}(1),\ldots,\mathbf{Y}(N)):=\prod_{k=1}^{K}
\widehat{E}[\widehat{\pi}_{Y_k}(Y_k)] & = & \prod_{k=1}^{K}\frac{1}{N}\sum_{n=1}^{N}\frac{1}{N}\sum_{m=1}^{N}\mathcal{K}_{2,h}\left
 ( \frac{Y_k(n)-Y_k(m)}{\sigma_{Y_k}}\right )\\
 & = & \frac{C_N^{2K}(2K)!}{N^{2K}} U_2 + \frac{b_N^{(2)}}{\sqrt{N}}\\
\end{eqnarray*}
where 
$$U_2(\mathbf{Y}(1),\ldots, \mathbf{Y}(N))=\frac{1}{N(N-1)\ldots (N-2K+1)}\sum_{S^{\neq}_2}\prod_{k=1}^{K}\mathcal{K}_{2,h}\left
 ( \frac{Y_k(i_k)-Y_k(j_{k})}{\sigma_{Y_k}}\right )$$
where ${S^{\neq}_2}$ is the subset of $S_2$ of 2K-dimensional elements whose all components are different from one another
and
$$b_N^{(2)}=\frac{1}{\sqrt{N}N^{(K-1)}}\sum_{S^=_2}\prod_{k=1}^{K}\mathcal{K}_{2,h}\left
 ( \frac{Y_k(i_k)-Y_k(j_k)}{\sigma_{Y_k}}\right )$$
and $S^=_2=S_2\setminus S^{\neq}_2$, that is, at least two components of each element are equal. 
\end{lem}

\begin{lem}
Let us define the set $S_3=\{(i_1,\ldots,i_{K},j)| \forall k, 1\leq k \leq K, 1 \leq i_k \leq N,  1 \leq j \leq N\}$.
Finally, the estimation of $E\left [\prod_{k=1}^{K}\pi_{Y_k}(Y_k)\right ]$, is expressed as:

\begin{eqnarray*}
U'_3(\mathbf{Y}(1),\ldots,\mathbf{Y}(N)):=\widehat{E}[\prod_{k=1}^{K}\widehat{\pi}_{Y_k}(Y_k)] & = & \frac{1}{N}\sum_{n=1}^{N}\prod_{k=1}^{K}\frac{1}{N}\sum_{m=1}^{N}\mathcal{K}_{2,h}\left
 ( \frac{Y_k(n)-Y_k(m)}{\sigma_{Y_k}}\right )\\
 & = & 2\frac{C_N^{K+1}(K+1)!}{N^{K+1}} U_3 + \frac{b_N^{(3)}}{\sqrt{N}}\\
\end{eqnarray*}
where
$$U_3(\mathbf{Y}(1),\ldots,\mathbf{Y}(N))=\frac{1}{N(N-1)\ldots (N-K+2)}\sum_{S^{\neq}_3}\prod_{k=1}^{K}\mathcal{K}_{2,h}\left
 ( \frac{Y_k(i_k)-Y_k(j)}{\sigma_{Y_k}}\right )$$
${S^{\neq}_3}$ is the subset of $S_3$ of (K+1)-dimensional elements whose all components are different from one another
and $$b_N^{(3)}=\frac{1}{\sqrt{N}N^{(K-1)}}\sum_{S^{=}_3}\prod_{k=1}^{K}\mathcal{K}_{2,h}\left
 ( \frac{Y_k(i_k)-Y_k(j)}{\sigma_{Y_k}}\right )$$
and $S^=_3=S_3\setminus S^{\neq}_3$, that is, at least two components of each element are equal. 
\end{lem}

\section{Proof of the computation of variance}\label{proof.lem.variance}

Using the decomposition of $\Qhatsigma$ in terms of U statistics, and following the computation of the variance given in \cite{hoeffding.1948.1}, the variance of $\Qhatsigma$ is computed.

\section{Proof of lemma \ref{lem.ind}}\label{proof.lem.ind}

Using the Chebychev inequality, we deduce :

$$1-P_{\text{H}_1}(\Qhat<q_{\alpha})=P_{\text{H}_1}(\Qhat>q_{\alpha})> 1-\frac{\text{var}(\Qhat)}{q_\alpha-Q}$$

Thanks to the asymptotic properties of the variance, the right hand side of the above expression tends to 1 when the sample size goes to infinity.

\end{document}